\input amstex
\documentstyle{amsppt}

\magnification=\magstephalf
\pagewidth{16true cm}
\pageheight{24.5true cm}
\frenchspacing
\parskip4pt

\topmatter

\title The manifold of finite rank projections 
in the space ${\Cal L} (H).$\endtitle

\author Jos\'e M. Isidro \dag\endauthor

\affil Universidad de Santiago de Compostela \endaffil

\address Facultad de Matem\'aticas, 
Universidad de Santiago,
15706 Santiago de Compostela, SPAIN \endaddress

\email jmisidro\@zmat.usc.es \endemail

\date  December, 1999\enddate

\thanks \dag Supported by Comisi\'on Hispano-H\'ungara 
de Cooperaci\'on Cient\'{\i}fica
y Tecnol\'ogica \endthanks


\subjclass 48 G 20,  72 H 51 \endsubjclass
 
\abstract 
Given a complex Hilbert space $H$ and the von Neumann algebra ${\Cal L}(H)$ 
of all bounded linear operators in $H$, 
we study the Grassmann manifold $M$ of all projections in ${\Cal L}(H)$ that 
have a fixed finite rank $r$. To do it we take the Jordan-Banach triple 
(or JB$^*$-triple) approach which allows us to define a natural Levi-Civita 
connection on $M$ by using algebraic tools. We identify the geodesics and  
the Riemann distance and establish some properties of $M$. 
\endabstract

\endtopmatter

\document 

\head{0} Introduction \endhead

In this paper we are concerned with the differential geometry of the 
infinite-dimensional Grassmann manifold 
$M$ of all projections in $Z\colon = 
{\Cal L}(H)$, the space of 
bounded linear operators $z\colon H\to H$ in a complex 
Hilbert space $H$. Grassmann manifolds are a classical object in 
Differential Geometry and in recent years several authors have considered 
them in the Banach space setting. Besides the Grassmann structure, a 
Riemann and a K\"ahler structure has sometimes been defined even in the 
infinite-dimensional setting. Let us recall some aspects of the history of the 
topic that are relevant for our purpose. 

The study of the manifold of minimal projections in a finite-dimensional simple 
formally real Jordan algebra was made by U. Hirzebruch in \cite{4}, who 
proved that such a manifold is a compact symmetric Riemann space of rank 1, and 
that every such a space arises in this way. Later on, Nomura in \cite{13, 14} 
established similar results for the manifold of fixed finite rank projections in a 
topologically simple real Jordan-Hilbert algebra. On the other hand, the 
Grassmann manifold $M$ of all projections in the space $Z\colon ={\Cal L}(H)$ 
of bounded linear operators has been discussed by Kaup in \cite{7} and \cite{10}. 
It is therefore reasonable to ask whether a Riemann structure  can always 
be defined in $M$ and how does it behave when it exists. It is known that $M$ 
has several connected components $M_r\subset M$ each of which consists of the 
projections in ${\Cal L}(H)$ that have a fixed rank $r$, $1\leq r
\leq \infty$. We prove that $M_r$ admits a Riemann structure if and only if 
$r<\infty$ establishing a distinction between the finite and the infinite 
dimensional cases. We then assume $r<\infty$ and proceed to discuss 
the behaviour of the Riemann manifold $M_r$, which looks very much like in 
the finite-dimensional case. One of the  novelties is that 
we take JB$^*$-triple approach instead of the Jordan-algebra approach of 
\cite{4} and \cite{13}. As noted in  \cite{1} and \cite{5}, within this 
context the algebraic structure of JB$^*$-triple acts as a substitute for the 
Jordan algebra structure and provides a {\sl local scalar product} known as 
the Levi form \cite{10}. Although ${\Cal L}(H)$ is not 
a Hilbert space, the JB$^*$-triple approach 
and the use of the Levi form allows us to define a torsion-free affine 
connection $\nabla$ 
on $M_r$ that is invariant under the group $\hbox{\rm Aut}^{\circ}(Z)$ of 
all surjective linear isometries of ${\Cal L}(H)$. We integrate the equation 
of the geodesics and define an $\hbox{\rm Aut}^{\circ}(Z)$-invariant 
Riemann metric on $M_r$ with respect to which $\nabla$ is a Levi-Civita 
connection. We prove that any two distinct points 
in $M_r$ can be joined by a geodesic which (except for the case of a pair of 
antipodal points) is uniquely determined and is a minimizing curve for 
the Riemann distance, that is also computed. We prove that $M_r$ is a symmetric 
manifold on which $\hbox{\rm Aut}^{\circ}(Z)$ acts transitively as a group 
of isometries.

\head{1} JB$^*$-triples and tripotents.\endhead 
For a complex Banach space $Z$, denote by ${\Cal L}(Z)$ the Banach algebra  
of all bounded linear operators on $Z$. A complex Banach space $Z$ with a 
continuous mapping $(a, b, c)  
\mapsto \{a b c\}$ from $Z\times Z\times Z$ to $Z$ is called a {\it JB*-triple} 
if the following conditions are satisfied for all $a, b, c, d \in Z$, where 
the operator $a\square b\in {\Cal L}(Z)$ is defined by $z\mapsto \{abz\}$ and 
$\lbrack\, , \, \rbrack$ is the commutator product:

\item{1} $\{abc\}$ is symmetric complex linear in $a, c$ and conjugate linear 
in $b$.
\item{2} $\lbrack a\square b , \, c\square d \rbrack = \{abc\}\square d -
c\square \{dab\}.$
\item{3} $a\square a$ is hermitian and has spectrum $\geq 0.$
\item{4} $\Vert \{aaa\}\Vert = \Vert a\Vert ^3$. 
\par
If a complex vector space $Z$ admits a JB*-triple structure, then the norm and 
the triple product determine each other.  
A {\it derivation} of a JB*-triple $Z$ is an element $\delta \in {\Cal L}(Z)$ such that 
$
\delta \{z z z \}= \{(\delta z)  z  z\} + \{z  (\delta z) z\} +
\{z z (\delta z)\} \; 
$ 
and an {\it automorphism} is a bijection $\phi \in {\Cal L}(Z)$ such 
that $\phi \{zz z\}= \{(\phi z) ( \phi z)(\phi z)\}$ for $z\in Z$. The latter
occurs if and only if $\phi$ is a surjective linear isometry of $Z$. 
The group $\hbox{\rm Aut}(Z)$ of automorphisms of
$Z$ is a real Banach-Lie 
group whose Banach-Lie algebra is the set of derivations of $Z$ 
. The connected component of the identity in $\hbox{\rm Aut}(Z)$ is denoted by 
$\hbox{\rm Aut}^{\circ}(Z)$. 
Two elements $x, y \in Z$ are {\sl orthogonal} if $x \square y = 0$.
An element $e\in Z$ is called a {\it tripotent} if $\{e e e \}=e$. The set 
$\hbox{\rm Tri}(Z)$ of tripotents  
is endowed with the induced topology of  $Z$. 
If $e\in \hbox{\rm Tri}(Z)$, then $e\square e\in {\Cal L}(Z)$ has the 
eigenvalues $0,\, {1\over 2}, \, 1$ and we have the topological direct sum 
decomposition
$$
Z=Z_1(e)\oplus Z_{1/2}(e)\oplus Z_0(e)
$$
called the {\it Peirce decomposition} of $Z$. Here $Z_k(e)$ is the $k$-
eigenspace and the {\sl Peirce projections} are 
$$
P_1(e) = Q^2(e), \qquad P_{1/2}(e)= 2(e\square e- Q^2(e)), \qquad 
P_0(e)= \hbox{\rm Id} -2e\square e+Q^2(e),
$$
where $Q(e)z = \{ e z e\}$ for $z\in Z$. 
We will use the {\sl Peirce rules}
$ \{Z_i (e)\, Z_j (e)\, Z_k (e) \} \subset Z_{i-j+k} (e) $
where $Z_l (e) = \{0\}$ for $l \neq 0,1/2,1$.
We note that $Z_1(e)$ is a complex unital JB*-algebra in 
the product $a\circ b\colon =
\{ a eb\}$ and involution $a^{\#} \colon =\{eae\}$. 
 Let
$$
A(e)\colon = \{ z\in Z_1(e)\;\colon \; z^{\#}=z\}.
$$
Then we have $Z_1(e)=A(e)\oplus iA(e)$. The {\sl Peirce spaces} of 
$Z$ with respect to a an orthogonal family of tripotents 
${\Cal E}=(e_i)_{i\in I}$ are defined by 
$$\align
Z_{ii}\colon &= Z_1(e_i) \\
Z_{ij}\colon &= Z_{1/2}(e_i)\cap Z_{1/2}(e_j), \quad i\neq j\\
Z_{i0}\colon &= Z_{0i}\colon = Z_{1/2}(e_i)\bigcap_{j\neq i}Z_0(e_j) \\
Z_{00}\colon &= \bigcap_{i\in I}Z_0(e_i)
\endalign
$$
The {\sl Peirce sum} $P({\Cal E})\colon = \bigoplus_{i,j\in I}Z_{ij}$ 
relative to the family ${\Cal E}$ is direct and we have $Z= P({\Cal E})$ 
whenever ${\Cal E}$ is a finite set. Every ${\Cal E}$-Peirce space is 
a JB$^*$-subtriple of $Z$ and the Peirce rules 
$$ \{ Z_{ij}Z_{jk}Z_{kl} \}\subset Z_{il}$$
hold for all $i,\,j,\,k,\,l\in I$. 

A tripotent $e$ in a JB$^*$-triple $Z$ is said to be {\sl minimal} if 
$ P_1(e)Z=\Bbb C e$, and we let $\hbox{\rm Min}(Z)$ be the set of them. Clearly $e=0$ 
lies in $\hbox{\rm Min}(Z)$ and is an isolated point there. 
If $e\in \hbox{\rm Min}(Z)$ and $e\neq 0$ then $\Vert e\Vert=1$ and by the Peirce multiplication 
rules we have $\{ euv\}\in Z_1(e)= \Bbb C e$ for all $u, v\in Z_{1/2}(e)$. 
Therefore we can define a sesquilinear form, called the {\sl Levi form}, 
$\langle \cdot , \cdot \rangle _e\colon Z_{1/2}(e)\times Z_{1/2}(e)\to \Bbb C$ by 
$$\{ e uv\} = \langle v , u \rangle_e \, e,\qquad u, v\in Z_{1/2}(e).
$$
It is known \cite{10} that $\langle \cdot , \cdot \rangle _e$ is positive 
definite hence it defines a scalar product in $Z_{1/2}(e)$ whose norm, called the 
{\sl Levi norm} and denoted by $\vert \cdot \vert_e$, satisfies  
$$\vert u\vert^2_e\leq  \Vert u\Vert ^2, \qquad u\in Z_{1/2}(e)$$ 
that is, we have the continuous inclusion $(Z_{1/2}(e), \, \Vert \cdot \Vert)
\hookrightarrow (Z_{1/2}(e), \, \vert \cdot \vert_e).$ To simplify the 
notation, we shall omit the subindex $e$ in both the Levi form and the Levi norm if no confusion 
is likely to occur. 

JB*-triples include C*-algebras and JB*-algebras.  
A C*-algebra is a JB*-triple with respect to the triple product 
$2\{abc\} \colon = 
(ab^*c+cb^*a).
$ 
Every JB*-algebra with Jordan product $(a, b)\mapsto 
a\circ b$ and involution $a\mapsto a^*$ is a JB*-triple with triple product
$\{abc\} = (a\circ b^*)\circ c -
(c\circ a)\circ b^* + (b^* \circ c)\circ a.
$ 

We refer to \cite{8,9,10,12} for the background of  
JB$^*$-triples theory.

\head{2} The manifold $M$ of minimal projections \endhead

Let $Z\colon = {\Cal L}(H)$, where $H$ is a complex Hilbert space, and 
let $M\subset \hbox{\rm Tri}(Z)$ denote the set of all 
projections in $Z$ 
endowed with its topology as subspace of $Z$. Fix any non zero  
projection $e_0\in M$ and denote by $M$ the connected component 
of $e_0$ in $M$. Then all elements in $M$ have the same 
rank as $e_0$ and 
$\hbox{\rm Aut}^{\circ }(Z)$ acts transitively on $M$ which is an 
$\hbox{\rm Aut}^{\circ }(Z)$-
invariant real analytic manifold whose tangent space at a point $e\in M$ is 
$$ T_eM= Z_{1/2}(e)_s ,$$
the selfadjoint part of the ${1\over 2}$-eigenspace of $e$. If we set 
$k_u\colon = 2(u\square e-e\square u)$, then by \cite{1, th. 3.3} a local 
chart of $M$ in a suitable neighbourhood $U$ of $0$ in $Z_{1/2}(e)_s$ is 
given by
$$u\mapsto f(u)\colon = \exp \, k_u(e) .$$
Let ${\frak D}(M)$ be the Lie algebra of all real analytic vector fields on $M$, 
and as in \cite{1}, define an affine connection $\nabla$ on $M$ by 
$$(\nabla _X\,Y)_e\colon = P_{1/2}(e) Y ^{\prime}_eX_e, \qquad e\in M, 
\qquad X, Y\in {\frak D}(M)\eqno(1).$$
Then $\nabla$ is a torsion-free $\hbox{\rm Aut}^{\circ}(Z)$-invariant 
affine connection on $M$. 
For each $e\in M$ and $u\in Z_{1/2}(e)_s$ we let $\gamma_{e,u}\colon \Bbb R \to M$ 
denote the curve $\gamma_{e,u}(t)\colon = \exp \, t k_u(e)$. Clearly we have 
$\gamma_{e,u}(0)=e$ and ${\dot \gamma}_{e,u}(0)=u\in T_eM$. By \cite{1, th. 2.7}, 
$\gamma_{e,u}$ is a $\nabla$-geodesic of $M$. 
Let us introduce a 
binary product in $Z$ by $x\circ y\colon = \{x e y\}$. Then $(Z, \circ )$ 
is a complex Jordan algebra where, as usual, $x^{(n)}$ denotes the $n$-th 
power of $x$ in $(Z, \circ )$ for $n\in \Bbb N$. 
For $u\in Z_{1/2}(e)$, the real Jordan subalgebra of $(Z, \circ )$ generated by 
the pair $(e, u)$ is denoted by $J[e,u]$ and we have $\gamma _{e, u}(\Bbb R 
)\subset J[e, u]$. 

To make a more detailed study of the manifold $M$,  
we shall assume that $e_0$ is {\it minimal}. In such a case 
$J[e,u]$ coincides with the closed real linear span 
of the set $\{ e, \,u,\, u^{(2)}\}$, in particular $\dim J[e,u]\leq 3$ and 
$$
\gamma _{e, u}(t)=(\cos ^2 t\theta) \,e + ({1\over 2\theta}\sin 2t\theta )\,u
+({1\over \theta ^2} \sin ^2 t\theta )\, u^{(2)}, \qquad t\in \Bbb R
\eqno(2)$$
for some  angle $0\leq \theta <{\pi \over 2}$. If $a,b$ are two distinct 
minimal projections and they are not orthogonal (that is, if the Peirce 
projection $P_1(a)b$ is invertible in the JB$^*$-algebra $Z_1(a)$) then 
there is an unique geodesic $\gamma_{a,u}(t)$ joining $a$ with $b$ in 
$M$. Moreover, due to the minimality of $e$ the tangent space $Z_{1/2}(e)
\approx \{e\}^{\perp}$ appears naturally endowed with the Levi form $\langle \cdot 
,\rangle_e$ and it turns out that the Levi norm $\vert \cdot\vert _e$ and 
the operator norm $\Vert \cdot\Vert$ are equivalent in $Z_{1/2}(e)$ (see 
\cite{6, th.5.1}). Thus $(Z_{1/2}(e), \, \vert\cdot \vert_e)$ is a Hilbert space and an  
$\hbox{\rm Aut}^{\circ}$-invariant Riemann structure can be defined in $M$ 
by
$$ g_e(X,Y)\colon = \langle X_e, Y_e\rangle_e, \qquad X,\, Y\in {\frak D}(M)
\eqno(3)
$$
where $V_e\in Z_{1/2}(e)$ denotes the value taken by the vector field $V$ 
at the point $e\in M$. By \cite{1} $g$ satisfies
$$
Xg(Y,Z) =g(\nabla_X Y,\,Z)+g(Y,\, \nabla_XZ), \qquad X,\,Y,\, Z\in 
{\frak D}(M) \eqno(4)
$$
Therefore $\nabla$ is the only Levi-Civita affine connection on $M$, and 
the geodesics are minimizing curves for the Riemann distance in $M$, 
which is given by the formula
$$
d( a, b) = \cos ^{-1} \big( \Vert P_1(a)b\Vert ^{1\over 2}\big)=\theta.
$$
$M$ is symmetric Riemann manifold on which $\hbox{\rm Aut}^{\circ}(Z)$ acts 
transitively as a group of isometries and there is a real analytic 
diffeomorphism of $M$ onto the projective space $\Bbb P (H)$ 
over $H$, endowed with the Fubini-Study metric. We refer to \cite{1,5,6,13} 
for proofs and background about these facts.

\head{3} The manifold of finite rank projections in ${\Cal L}(H).$ \endhead 

In what follows we let $M$ and $M_r$ be the set of all projections in $Z$ 
and the set of all projections that have a fixed finite rank $r$, respectively. 
If $a\in M_r$ then a {\it frame} for $a$ is any family $(a_1, \cdots ,a_r)$ 
of pairwise orthognal minimal projections in $Z$ such that $a=\Sigma a_k$. Note 
that then the $a_k$ have the form $a_k=(\cdot ,\alpha_k )\alpha_k$ where 
$(\alpha_k)$ is an orthonormal family of vectors in the range $a(H)$. 
\proclaim{3.1 Proposition}
For every projection $a\in M$ the following conditions are equivalent:
\roster
\item The rank of $a$ is finite. 
\item The Banach space $Z_{1/2}(a)$ is linearly homeomorphic to a Hilbert space.
\endroster
\endproclaim
\demo{Proof} Let us choose an orthonormal basis $(\alpha_\imath )_{\imath\in I}$ 
in the range $a(H)\subset H$ of $a$. Then $a_\imath \colon = (\cdot , 
\alpha_\imath )\alpha_\imath , \, \imath\in I,$ is a family of pairwise orthogonal 
minimal projections that satisfy  
$$a=\Sigma_{\imath \in I}a_\imath \qquad \hbox{\rm strong operator 
convergence in $Z$}
\eqno(5)$$
The space $Z_{1/2}(a)_s$ 
consists of the operators $u\in Z$ such that $2\{aau\}=u$ and using (5) it 
is easy to check that $u$ can be represented in the form 
$$u=\Sigma_{\imath \in I}(\cdot , \xi_\imath )\alpha_\imath +(\cdot , 
\alpha_\imath )\xi_\imath \qquad \hbox{\rm strong operator convergence in $Z$}
$$
where $\xi_\imath \colon = 
u(\alpha_\imath)$ are vectors in $H$ that satisfy $\xi_\imath \in a(H)^{\perp}$. 
By (4) each $u\in Z_{1/2}(a)_s$ is determined by the family $(\xi_\imath )_
{\imath \in I}$. To simplify the notation, set $K\colon = a(H)^{\perp}$ and 
$L\colon = \ell _\infty (I, K)$ for the Banach space of the families 
$(\xi_\imath )_{\imath \in I}\subset K$ with the norm of the supremun 
$\Vert (\xi_\imath )\Vert \colon =\sup _{\imath \in I}\Vert \xi_\imath \Vert$. 
Then the mapping
$$ L\to Z_{1/2}(a)_s, \qquad (\xi_\imath )\mapsto u_\xi \colon = 
\Sigma _{\imath \in I}\,[ \,(\cdot , \alpha _\imath )\xi_\imath 
+(\cdot , \xi_\imath )\alpha_\imath \,]$$
is a continuous {\it real linear} vector space isomorphism, hence a homeomorphism . 
Thus if the operator norm in $Z_{1/2}(a)_s$ is equivalent to a Hilbert space norm 
the same must occur with $\ell _\infty (I, K)$, hence $I$ must be a finite set 
which means that $a= \Sigma a_\imath$ has finite rank.  The converse is 
easy. \qed
\enddemo
\proclaim{3.2 Lemma}
Let $a,\,b\in M_r$ with $a=\Sigma a_k$ where the $(a_k)$ is a frame 
for $a$, and let $Q(a_k)b= \lambda_ka_k$, 
$(k=1,\cdots ,r)$. If $P_1(a)b$ is invertible in the JB$^*$-algebtra
$Z_1(a)$, then $\lambda_k\neq 0$ for all $k$. The set of all elements 
$b\in M_r$ for which $P_1(a)b$ is invertible in $Z_1(a)$ is dense in $M_r$.
\endproclaim
\demo{Proof} Suppose that $a_k=(\cdot , \alpha_k)\alpha_k$ and 
$b_j=(\cdot , \beta_j)\beta_j$ are frames for $a$ and $b$ respectively. 
Then for each fixed $k$ we have
$$Q(a_k)b= \{a_kba_k\}= \big( \Sigma_j \vert (\alpha_k , \beta_j )\vert ^2
\big) a_k= \lambda _k a_k
$$
where $\lambda _k\geq 0$. Moreover $\lambda_k=0$ if and only if 
$\alpha_k\in \{\beta_1 , \cdots , \beta_r\}^{\perp}$ which is equivalent to 
$a_k\perp b$. But in such a case $\hbox{\rm range} (a_k)\subset 
\hbox{\rm ker} \{a_kba_k\} = \hbox{\rm ker} P_1(a)b$ which contradicts the 
invertibility of $P_1(a)b$. To simplify the notation set $K\colon = a(H)\subset 
H$ and note that $\hbox{\rm dim}\,K=\hbox{\rm rank}\, a=r<\infty$ The 
operators in $Z_1(a)= aZa$ can be viewed as operators in ${\Cal L}(K)$, 
therefore the {\it determinant} function is defined in $Z_1(a)$ and an element 
$z\in Z_1(a)$ is invertible if and only if $\hbox{\rm det}\, (z)\neq 0$. Thus 
the set of the operators $b\in Z$ for which $P_1(a)b$ is invertible in 
$Z_1(a)$ is an open dense subset of $M_r$. \qed 
\enddemo  
\proclaim{3.3 Lemma} 
If $a$, $p$ and $q$ are projections in $M_r$ and 
$P_{1/2}(a)p=P_{1/2}(a)q$, then $p=q$.
\endproclaim
\demo{Proof} Take frames for $a,\,p,\,q$, compute $P_{1/2}(a)p= 
2(D(a\square a)- Q(a)^2)p$ and 
proceed similarly with $q$. An elementary exercise of linear algebra yields 
range (p)=range (q), hence $p=q$.\qed 
\enddemo 
Let $a\in M_r$ and choose any frame 
$(a_1, a_2, \cdots , a_r)$ for $a$. 
As above $Z_{1/2}(a)_s$ consists of the operators 
$u= \Sigma (\cdot , \xi_k )\alpha_k+(\cdot , 
\alpha_k)\xi_k$
where $\xi_k\colon =u(\alpha_k)$ are vectors in $H$ 
that satisfy $\xi_k\in a(H)^{\perp}.$ 
Write $u_k\colon = (\cdot , \xi_k )\alpha_k+(\cdot , \alpha_k)
\xi_k$. Then we have 
$u= \Sigma u_k $
where the $u_k$ are selfadjoint operators in $Z={\Cal L}(H)$ 
(in fact $u_k\in Z_{1/2}(a_k)_s$) that satisfy
$$ 
u_j\square a_k=a_k\square u_j=0, \qquad  j\neq k , \qquad (j, k = 1,2,
\cdots , r)\eqno(6). 
$$

The above properties of the $a_k,\,u_k$ hold whatever is the frame $(a_1, 
a_2,\cdots , a_r)$. There are many families in those conditions and 
we are going to prove that, by making an appropriate choice of the $a_k$ 
(a choice in which the tangent vector $u\in Z_{1/2}(a)$ is also involved) 
we can additionally have 
$$u_k\square u_j=u_j\square u_k=0, \qquad j\neq k , \qquad (j,k = 
1,2,\cdots , r) \eqno(7)$$
This will simplify 
considerably the calculations in the sequel. 
We need some material. 
\proclaim{3.4 Lemma} With the above notation the set of minimal tripotents 
in $Z_{1/2}(a)$ is 
$$
\{ \,(\cdot , \alpha )\xi+(\cdot , \xi )\alpha \;\colon \;
\alpha \in a(H), \; \xi \in a(H)^{\perp}, \; \Vert \alpha \Vert = 
1=\Vert \xi\Vert \,\}$$
\endproclaim 
\demo{Proof}
Let $x\in Z$ be of the form $x=(\cdot , \alpha )\xi+(\cdot , \xi )\alpha$ 
where $\alpha , \xi\in H$ satisfy the above conditions. It is a matter of 
routine calculation to see that then $2\{aax\}=x$ hence $x\in Z_{1/2}(a)_s$.  
Moreover $\{xxx\}=x$ so that $x$ is a tripotent and we can easily see that 
$\{ x Z_{1/2}(a)x\}\subset \Bbb C x$ which proves the minimality of $x$ in 
$Z_{1/2}(a)$. The converse is similar. \qed
\enddemo
The following result should be compared to \cite{14, prop. 3.4} 
\proclaim{3.5 Lemma} Two minimal tripotents $x=(\cdot , \alpha )\xi+
(\cdot , \xi )\alpha$ and $y=(\cdot , \beta )\eta+
(\cdot , \eta )\beta$ in $Z_{1/2}(a)_s$ are orthogonal if and only if 
$\alpha \perp \beta$ and $\xi \perp \eta$. In particular $Z_{1/2}(a)$ has 
rank $r$ for all $a\in M$
\endproclaim  
\demo{Proof}
By \cite{2, p. 18} $x$ and $y$ are orthogonal if and only if the conditions 
$xy^*=0=y^*x$ hold. Now it is elementary to complete the proof of the 
first statement. For the second part, let $(u_\imath )_{\imath \in I}$ 
be a family of pairwise minimal orthogonal tripotents in $Z_{1/2}(a)$. 
Then $u_{\imath}=(\cdot , \alpha_{\imath} )\xi _{\imath}+ (\cdot , 
\xi _{\imath} )\alpha_{\imath }$ where $(\alpha_{\imath })\subset a(H)$ and 
$\xi _{\imath}\subset a(H)^{\perp}$ are orthonormal families of vectors 
in $H$. In particular $a_{\imath }\colon = (\cdot , \alpha _{\imath})
\alpha _{\imath}$ is a family of pairwise orthogonal projections with 
$\Sigma a_{\imath}\leq a$. Since rank(a)=r, we have cardinal$\,(I) \leq r$. 
The converse is easy. \qed 
\enddemo
Let $a\in M$ be a fixed projection and take any tangent vector $u\in 
Z_{1/2}(a)_s$ to $M$ at $a$. 
By lemma 3.2 $Z_{1/2}(a)$ has finite rank, hence  
\cite{9, cor 4.5} $u$ has a spectral decomposition 
in the JB$^*$-triple $Z_{1/2}(a)$ of the form 
$$u= \rho_1 u_1+\cdots +\rho_su_s, \qquad 0\leq\rho_1\leq \cdots 
\leq\rho_s=\Vert u\Vert , \qquad 1\leq s\leq r \eqno(8)$$
where the $u_k$ are pairwise {\it orthogonal minimal tripotents} in $Z_{1/2}(a)$.  
Therefore 
$$
\align
u_k=(\cdot &, \alpha_k )\xi_k+(\cdot , \xi_k )\alpha_k , 
\;\;\alpha_k\in a(H), \;\; \xi_k\in a(H)^{\perp}, \\
 \Vert \alpha_k\Vert &= 1=
\Vert \xi_k\Vert, \;\; \alpha _j\perp \alpha_k , \;\; \xi_j\perp \xi_k, \;\;
j \neq k 
\endalign
$$ 
Then $a_k\colon = (\cdot , \alpha_k)\alpha_k$ are pairwise orthogonal 
minimal projections in $Z$ and $\Sigma a_k\leq a$. In case $s<r$, which 
occurs if some of the $\rho_k=0$, we pick additional minimal orthogonal 
projections $a_{s+1}, \cdots , a_r$ so as to have $a=\Sigma a_k$. For the family 
$(a_1,\cdots ,a_r)$ so constructed, called a {\it frame} 
associated to the pair $(a,\,u)$, both properties (6) and (7) hold. Remark 
that this frame needs not be unique, it depends on $a$ and on $u$ as well, and it is 
invariant under the group $\hbox{\rm Aut}^{\circ} (Z)$. In fact some more 
properties are valid now. 

In accordance with section \S 1, each pair $(a_k, u_k)$ gives rise to a real Jordan algebra 
$J_k \colon =J[a_k, u_k]$ with the product $x\circ_{k} y \colon = \{ xa_ky\}$. 
We have $\hbox{\rm dim} (J_k)=3$ and $\{ a_k, u_k, u_k^{(2)}\}$ is a basis of $J_k$. 
Moreover, $J_k$ is invariant under the operator 
$g_k\colon = 2( a_k\square u_k-u_k\square a_k)$ where triple products are 
computed in $Z={\Cal L}(H)$. In case $s<\hbox{\rm rank}(a)$ we set $J_n 
\colon =\Bbb R a_n$ as real Jordan algebras. 

\proclaim{3.6 Lemma} The Jordan algebras $J_k$ and $J_l$ with $k\neq l$,
 $(k, l=1, \cdots , r)$ are orthogonal in the JB$^*$-triple sense in $Z$, 
that is $\{J_kJ_lZ\}=0$.
\endproclaim  
\demo{Proof} 
For $n\in\{ k, l\}\subset \{1, \cdots , s\}$ with $k\neq l$, let $z_n$ be 
any element in the basis $\{a_n, u_n, u_n^{(2)}\}$ of $J_n$. Clearly it 
suffices to show that $z_kz_l=0=z_lz_k$. As an example, we shall prove that 
$u_k^{(2)}u_l^{(2)}=0$. 
It is a routine to check that 
$u_ku_l=0$. Then 
$$u_k^{(2)}u_l^{(2)}= \{ u_ka_ku_k\}\, \{u_la_lu_l\}= (u_ka_ku_k)\,
(u_la_lu_l) = u_ka_k(u_ku_l)a_lu_l=0$$
as we wanted to see. \qed
\enddemo
Consider now the vector space direct sum $J\colon = \bigoplus_1^r J_k$, and 
define a product $z\circ w\colon = \{zaw\}$ in $J$ by 
$$z\circ w\colon =\{zaw\}= {1\over 2}(zaw+waz)= {1\over 2}\Sigma_1^r(z_ka_k
w_k+w_ka_kz_k)= \Sigma_1^r z_k\circ_k w_k$$ 
where $z_k, \,w_k$ are respectively the $J_k$-component of $z$ and $w$. It 
is now clear that $J$ is a real Jordan algebra, that the product in $J$ induces 
in each $J_k$ its own product $z\circ_k w=\{za_k w\}$ and that the $J_k$ 
are orthogonal as Jordan subalgebras of $J$. It is also clear that $J$ 
coincides with the closed real linear span of the set $\bigcup_1^r 
\{a_k,\, u_k,\,u_k^{(2)}\}$, in particular $\hbox{\rm dim} J\leq 3r<\infty$. Finally 
$J[a,u]\subseteq J$ and we conjecture that the equality holds (see \cite{14, 
prop. 3.5 \& th. 3.6}.

\head{4} Geodesics and the exponential mapping. \endhead 

Consider $M_r$ endowed with the affine connection $\nabla$ given by (1). 
To discuss its geodesics, let us define an operator $g\in Z={\Cal L}(H)$ by 
$$ g\colon =g_{a,u}\colon = 2(u\square a -a\square u)= 2\Sigma \rho_k 
(u_k\square a_k-a_k\square u_k)= \Sigma \rho_k g_{a_k, u_k}$$
where $u=\Sigma \rho_ku_k$ is the spectral decomposition of $u\in 
Z_{1/2}(a)$, the $a_k$ is any frame associated to the pair $(a, u)$ and 
$g_k\colon =g_{a_k, u_k}$ is defined in a obvious manner. If the spectral 
decomposition of $u$ (see (8)) has $s<r$ non zero summands then we define 
$g_n\colon =0$ for $n=s+1,\cdots ,r$. 
Then $g_k$ is a commutative family of operators in $Z$, more precisely we have 
$g_k(J_l)=\{0\}$, $g_k\, g_l=g_l\, g_k=0$ for all $k\neq l$, $(k, l=1,
\cdots , r)$ and $g$ leaves invariant all the spaces $J$ and $J_k$. Thus
$$
\gamma_{a,u}(t)\colon =\exp \, t g(a) =\Sigma \exp \, t g_k (a_k), 
\qquad t\in \Bbb R
$$ 
By section \S1 this curve is a geodesic in $M_r$ and $\gamma_{a,u}(\Bbb R)\subset
J[a,u]\subset J$.  
We can collect now the above discussion in the following 
statement (see \cite{14, prop. 5.1 \& 5.4}
\proclaim{4.1 Theorem} 
Suppose that we are given a point $a\in M_r$ 
and a tangent vector $u\in Z_{1/2}(a)_s$ to $M_r$ at $a$. 
Then the geodesic 
of $M_r$ that passes through $a$ with velocity $u$ is the curve 
$$\gamma _{a, u}(t)=\Sigma \gamma_{a_k, u_k}(t), \qquad t\in \Bbb R, $$
where $\gamma_k\colon =\gamma_{a_k,u_k}$ is given by 
$$\gamma_k(t)\colon =\gamma_{a_k, u_k}(t) = 
(\cos ^2 \theta_k t)\;a_k \;+({1\over 2\theta_k}\sin 2\theta_k t)\;u_k
+({1\over \theta_k^2} \sin ^2 \theta_kt)\; u_k ^{(2)} \eqno(G)$$
Here $u=\Sigma \rho_ku_k$ is the spectral decomposition of $u$ in 
$Z_{1/2}(a)$, the $a_k$ form a frame associated to the pair $(a,\,u)$ and 
the numbers $\theta_k$ are given by $\cos ^2\theta_k \colon = 
\rho_k$ with $0\leq \theta_k<{\pi \over 2}$. 
\endproclaim
Now we are in a position to define the exponential mapping. Suppose the 
tangent vector $u$ lies in the unit ball $B_1(a)\subset Z_{1/2}(a)$, i.e.  
$\Vert u\Vert <1$. For $t=1$ 
the expression (G) yields 
$$\gamma(1) =\Sigma 
(\cos ^2 \theta_k )\;a_k \;+\Sigma ({1\over 2\theta_k}
\sin 2\theta_k )\;u_k
+\Sigma ({1\over \theta_k^2} \sin ^2 \theta_k)\; u_k ^{(2)} \eqno(E)$$
and a real analytic mapping form the unit ball $B_1(0)\subset Z_{1/2}(a)$ 
to the manifold $M$ can be defined by 
$$\hbox{\rm Exp}_a(u)\colon = \gamma_{a,u}(1)$$
An inspection of (E) yields that the Peirce decomposition of $\gamma_{a,u}(1)$ 
relative to $a$ is  
$$\align 
P_1(a) \gamma_{a,u}(1)= \Sigma (\cos ^2\theta_k)a_k, &\qquad 
P_{1/2}(a)\gamma_{a,u}(1)=\Sigma ({1\over 2\theta_k}
\sin ^2 \theta_k)\; u_k \\
P_0(a) \gamma_{a,u}(1)&=
\Sigma ({1\over \theta_k^2} \sin ^2 \theta_k)\; u_k ^{(2)}
\endalign
$$
Remark that $0<\cos ^2\theta_k\leq 1$, hence in particular $P_1(a)\gamma_{a,u}(1)$ 
lies in the set of all ${\Cal N}_a$ of all 
invertible elements in the JB$^*$-algebra $Z_1(a)$. Clearly ${\Cal N}_a$ is 
an open neighbourhood of $a$ in $Z_1(a)$. 
Remark also that  
$0\leq {1\over 2\theta_k}\sin ^2 \theta_k=\rho_k \leq \Vert u\Vert <1$, hence 
$\Sigma ({1\over 2\theta_k}\sin ^2 \theta_k)\; u_k $ is the spectral 
decomposition of $P_{1/2}(a)\gamma_{a,u}(1)$ in $Z_{1/2}(a)$. 
Thus $\hbox{\rm Exp}_aB_1(a)\subset {\Cal N}_a\subset M$. We refer to 
$\hbox{\rm Exp}_a$ as the {\it exponential} mapping.

\head{5} Geodesics connecting two given points. The logaritm mapping. \endhead

Now we discuss the possibility of joining two given projections $a$ and $b$ 
such that $P_1(a)b$ is invertible in the Jordan algebra $Z_1(a)$, by means of 
a geodesic in $M$. The remarks in the precedent section show how to proceed. 
First we compute the spectral decomposition of $u\colon = P_{1/2}(a)b$ in 
the JB$^*$-triple $Z_{1/2}(a)$. Assume it to be 
$$ u= 
P_{1/2}(a)b= \Sigma \rho_ku_k, \qquad 0\leq\rho_1\leq \cdots \leq \rho_ r = 
\Vert u\Vert<1  , \qquad 1\leq k\leq r 
$$
where the $u_k$ are pairwise orthogonal minimal tripotents in $Z_{1/2}(a)$. 
Hence 
By lemma 3.4 the $u_k$ have the form $u_k= (\cdot , \alpha _k)
\xi_k+(\cdot , \xi_k)\alpha_k$ for some orthonormal families of vectors 
$(\alpha_k )\subset a(H)$ and $(\xi_k )\subset a(H)^{\perp}$. By lemma 3.2
$Q(a_k)b=\{a_kba_k\}= \lambda_k$ 
where $\lambda_k\neq 0$ since 
$P_1(a)b$ is invertible in $Z_1(a)$. Also $\vert \lambda_k\vert =
\Vert \{a_kba_k\}\Vert \leq 1$. Thus $0<\lambda_k\leq 1$ and a unique 
angle $0\leq \theta_k<{\pi \over 2}$ is determined by $\cos ^2 \theta_k = 
\lambda _k$. In this way we have got all the elements appearing in $(E)$. 
Let us define $\tilde \gamma (t) \colon = \Sigma \tilde \gamma _k(t)$ for 
$t\in \Bbb R$ where 
$$
\tilde \gamma _k(t) \colon = (\cos ^2 t \theta_k)\; a_k  + ({1\over 2\theta_k}
\sin 2t\theta_k )\;u_k+ ({1\over \theta_k^2} \sin ^2 t\theta_k )\; u_k^{(2)}
$$
By section \S 1, each $\tilde\gamma _k(t)$ is a geodesic in the manifold $M_1$ 
of all rank 1 projections. By the previous discussion $\tilde\gamma _j(t)$ and 
$\tilde\gamma_k(t)$ 
are orthogonal whenever $j\neq k$, $t\in \Bbb R$, hence 
$
\tilde\gamma (t)\colon = \Sigma \tilde\gamma _k(t), 
\; t\in \Bbb R ,
$
is a curve in the manifold $M$ of projections of rank r. Clearly 
$\tilde\gamma (0)=\Sigma \tilde\gamma_k(0)=\Sigma a_k=a$ and we shall now show that 
$\tilde b\colon =\gamma (1)$ coincides with $b$. As above 
$P_{1/2}(a) =\tilde b =\Sigma ({1\over2}
\sin 2\theta _k) \,u_k =\Sigma \rho_k u_k$ is the spectral decomposition of 
$P_{1/2}(a) \tilde b $ in $Z_{1/2}(a)$, which by construction is 
the spectral decomposition of $P_{1/2}(a)b$. 
Hence by lemma 3.3, $\tilde b =\tilde\gamma (1)=b$. This gives a geodesic $\gamma (t)$ 
that connects $a$ with $b$ in the manifold $M_r$ and passes through the point 
$a$ with the velocity $u\colon =P_{1/2}(a)b$. It is uniquely 
determined by the data $a, b$ and the property $\gamma_{a,u} (1)=b$. 

Now we are in a position to define the logaritm mapping. Fix a point 
$a\in M$ and let ${\Cal N}_a\subset M$ be the set of all projections $b\in M$ 
such that $P_1(a)b$ is invertible in the JB$^*$-algebra $Z_1(a)$. 
Define a mapping $\hbox{\rm Log}_a$ from ${\Cal N}_a\subset M$ to the 
unit ball $B_1(a)\subset  Z_{1/2}(a)$ by declaring $\hbox{\rm Log}_a(b)$ 
to be the velocity at $t=0$ of the unique geodesic $\gamma_{a,u}(t)$ that joins 
$a$ with $b$ in $M$ and $\gamma_{a,u}(1)=b$, in other words $\hbox{\rm Log}_a (b)
\colon = P_{1/2}(a)b$. We refer to $\hbox{\rm Log}_a$ as the 
{\it logaritm} mapping. Clearly $\hbox{\rm Log}_a$ and $\hbox{\rm Exp}_a$ are 
real analytic inverse mappings. In particular, the family $\{ ({\Cal N}_a, 
\hbox{\rm Log}_a)\, \colon \, a\in M\}$ is an atlas of $M$. We remark 
the fact that 
$\gamma_{a,u}[0,1]\subset {\Cal N}_a$ for all $u\in B_1(a)$ which shall be 
needed later on to 
apply the Gauss lemma \cite{11, 1.9} and  
summarize the above 
discussion in the statement (see \cite{14, th. 5.7 \& prop. 5.8})
\proclaim{5.1 Theorem} Let $a$ and $b$ be two given projections 
in $M_r$ and assume that $P_1(a)b$ 
is invertible in the Jordan algebra $Z_1(a)$. Then there is exactly 
one geodesic $\gamma_{a,u}(t)$ that joins $a$ with $b$ in $M$ and 
$\gamma_{a,u}(1)=b$. 
\endproclaim

\head{6} The Riemann structure on $M$.\endhead 
Let $a\in M_r$ and choose any frame $(a_k)$ for $a$. By section \S1 we 
have vector space direct sum decomposition 
$$
Z_{1/2}(a)= \bigoplus_1^rZ_{1/2}(a_k) \eqno(9)
$$ 
which suggests to define a scalar product in $Z_{1/2}(a)$ by 
$$
\langle u,\,v\rangle \colon ={1\over \sqrt r} \Sigma \langle u_k,\, v_k\rangle _{a_k} 
\eqno(10)$$
where $\langle \cdot , \cdot \rangle_{a_k}$ stands for the Levi form on 
$Z_{1/2}(a_k)$. First we prove 
\proclaim{6.1 Lemma}
With the above notation, $\hbox{\rm (9)}$ defines an $\hbox{\rm Aut}^{\circ}$-invariant 
scalar product on $Z_{1/2}(a)$ that does not depend of the frame 
$a=\Sigma_k$ and converts $Z_{1/2}(a)$ into a Hilbert space. 
\endproclaim
\demo{Proof}
Let $\Sigma a_k$ and $\Sigma a^{\prime}_k$ denote two frames for $a$ 
where $a_k=(\cdot , \alpha_k)\alpha_k$ and $a^{\prime}_k=(\cdot , 
\alpha_k^{\prime})\alpha_k^{\prime}$ for some orthonormal families $(\alpha_k),\,
(\alpha_k^{\prime})\subset a(H)$. Extend them to two orthonormal 
basis of $H$ and let $u\in {\Cal L}(H)$ be the unitary 
operator that exchanges these bases. Then $u$ induces an isometry $U\in \hbox{
\rm Aut}^{\circ}(Z)$ by $Uz= uzu^{-1}$ that satisfies $Ua_k^{\prime}=a_k$. The 
invariance of the Levi form together with (10) yields part of 
the result. The remainder is trivial. \qed 
\enddemo
A Riemann structure can now be defined in $M_r$ in the following way. Let $X,\, Y
\in {\frak D}(M)$ vector fields on $M_r$, and for $a\in M_r$ take any frame 
$a=\Sigma a_k$. Then (9) gives representation $X=\Sigma X_k, \, 
Y=\Sigma Y_k$ with $X_k,\, Y_k\in Z_{1/2}(a_k)$ and we set 
$$
g_a(X,\,Y)\colon =\langle X,\, Y\rangle ={1\over \sqrt r}\Sigma \langle X_k,
Y_k\rangle_{a_k}={1\over \sqrt r}\Sigma g_{a_k}(X_k,\, Y_k)
$$
This is a well defined $\hbox{\rm Aut}^{\circ}$-invariant Riemann structure 
on $M_r$. By section \S1 each $g_{a_k}$ has property (4) and a routine argument 
gives the same property for $g$. Thus $g$ is the only Levi-Civita connection 
in $M_r$ and we can apply the Gauss lemma \cite{11, 1.9} to conclude that 
the $\nabla$-geodesics are minimizing curves for the Riemann distance. 

Recall that for a tripotent $a\in Z$, the mapping 
$\sigma_a\colon x_1+x_{1/2}+x_0\mapsto x_1-x_{1/2}+x_0$, where $x\in Z$ and 
$x_1+x_{1/2}+x_0$ is the Peirce decomposition of $x$ with respect to $a$, called 
the Peirce symmetry of $Z$ with center $a$, is an involutory automorphism of $Z$ 
that induces an isometric symmetry of $M_r$ (see \cite{6, th. 5.1}). 
We let $\hbox{\rm Isom} M_r$ and ${\frak S}$ denote the group of all 
isometries of the Riemann manifold $M_r$ and the subgroup generated by the set 
$S\colon =\{\sigma_a \,\colon \, a\in M_r\}$, respectively. 
\proclaim{6.2 Proposition}  
With the above notation, $M_r$ is a symmetric Riemann manifold in which the 
group $\frak S$ acts transitively.
\endproclaim
\demo{Proof}
Let $a,\,b\in M_r$ be such that $b\in {\Cal N}_a$. Then $a$ and $b$ can be 
joined in $M_r$ by a unique geodesic with $\gamma (0)=a, \,\gamma (1)=b$. 
If $c\colon =\gamma ({1\over 2})$, then $\sigma_c$ is a symmetry of $M_r$ 
such that $\sigma_c(a)=b$. Thus the set $S$ is transitive in ${\Cal N}_a$ 
and $S$ is locally transitive in $M_r$. Consider now the case $b\notin {\Cal N}_a$. 
Since $M_r$ is pathwise connected, we can join $a$ with $b$ by a curve $\Gamma$ 
in $M_r$ and by a standard compactness argument there exists a finite set 
$\{ b_0,\cdots ,b_s\}\subset \Gamma [0,1]$ such that $b_0=a, \, b_s=b$ 
and $b_{k+1}\in {\Cal N}_{b_k}$ for $k=1,\cdots , s$. An application of the 
above argument to each pair of consecutive points gives the result. \qed 
\enddemo 
We now compute the Riemann distance in $M_r$. Consider first the case of two 
points $a, \,b\in M_r$ with $b\in {\Cal N}_a$. Let $\gamma_{a,u}(t)$ be the 
unique geodesic that joins $a$ with $b$ in $M_r$ and satisfies $b=\gamma_{a,u}(1)$. 
Since $\hbox{Aut}^{\circ} (Z)$ is transitive in ${\Cal N}_a$ and the Levi 
norm is $\hbox{Aut}^{\circ} (Z)$-invariant, we have
$$
\vert {\dot \gamma}_{au}(t)\vert_{\gamma_{au}(t)}= 
\vert {\dot \gamma}_{au}(0)\vert_{\gamma_{au}(0)} =
\vert u\vert_a
$$
On the other hands, since the Levi norm in $Z_{1/2}(a)$ is the direct 
hilbertian sum of the Levi norms in the $Z_{1/2}(a_k)$, we have by section \S1 
$$
\vert u\vert_a^2= {1\over r} \Sigma \vert u_k\vert_{a_k}^2= 
{1\over r }\Sigma \theta_k^2 \eqno(D)
$$
where $u=\Sigma \rho_ku_k$ is the spectral decomposition of $u$ in $Z_{1/2}(a)$, 
$(a_k)$ is the frame associated to the pair $(a,u)$ and $\cos ^2 \theta_k= 
\rho_k$. Therefore
$$
d(a,b)= \int_0^1 \vert {\dot \gamma}_{au}(t)\vert_{\gamma_{au}(t)}\;dt = 
\int_0^1 \vert u\vert_a\, dt= \vert u\vert_a= {1\over \sqrt r} 
\big( \Sigma \theta_k^2 \big) ^{1/2}
$$
Consider now the case $b\notin {\Cal N}_a$. By lemma 3.2 we can take a sequence 
$(b_n)_{n\in \Bbb N}$ in ${\Cal N}_a$ such that $b=\lim_{n\to \infty} b_n$.  
since (D) holds for all $b_n$ and the Riemann distance is continuous, we get 
the validity (D) for all $a,b\in M_r$. \qed 

Note that expression (D) is a generalization of the classical formula for the 
Fubini-Study metric in the projective space $\Bbb P (H)$.

\enddocument